\documentclass[11pt,leqno]{amsart}
\usepackage{amsmath,amssymb,amsthm}
\input{diagrams}
\begin{document}
\title[New components of the Moduli space]{New components of the Moduli space of minimal surfaces in $4$-dimensional flat tori\footnote{2000 Mathematics Subject Classification 49Q05, 53A10}}
\author{Toshihiro Shoda}
\address{Department of Mathematics, Tokyo Institute of Technology, 2-12-1 Oh-okayama, Meguro-ku, Tokyo 152-8551 Japan}
\email{tshoda@math.titech.ac.jp}
\maketitle
\begin{abstract}
In this paper, we consider new components of a key space of a Moduli space of minimal surfaces in flat $4$-tori and calculate their dimensions. Moreover, we construct an example of minimal surfaces in $4$-tori and obtain an element of the Moduli. In the process of the construction, we give an example of minimal surfaces with good property in a $3$-torus distinct from classical examples.
\end{abstract}

\newtheorem{defi}{Definition}[section]
\newtheorem{thm}{Theorem}[section]
\newtheorem{pro}{Proposition}[section]
\newtheorem{cor}{Corollary}[section]
\newtheorem{lem}{Lemma}[section]
\newtheorem{rem}{Remark}[section]
\newtheorem{main}{Main Theorem}

\section{Introduction}

Our objects are compact minimal surfaces of genus $g$ in $n$-dimensional flat tori, and the main studies are the following two questions:
 (i) To investigate Moduli spaces of minimal surfaces in flat tori,
 (ii) To give an example of minimal surfaces in flat tori.
 We consider these problems in the case $n=4$. Then, we find new components of the key space of the Moduli and calculate their dimensions. Moreover, we construct an example in the low genus case and obtain an element of the Moduli. The construction is outlined as follows. We first give a minimal surface in a flat $3$-torus. By a suitable deformation of the surface, we construct a minimal surface in a flat $4$-torus. We emphasize that we obtain a minimal surface with good property in a flat $3$-torus distinct from classical examples.

Now we explain the backgrounds and state our main results. First, we explain about (i). In a previous paper, we constructed a family of examples of minimal surfaces in flat $4$-tori \cite{S}. This family of minimal surfaces motivates us to study a set of minimal surfaces, that is, the Moduli space of minimal surfaces. The Moduli space $\mathcal{M}^n_g(\textbf{Q})/\mathcal{O}(n)$ of minimal surfaces in flat tori has been studied \cite{A-Pi}, \cite{Pi}. The key space is $\mathcal{M}^n_g$ which is defined as a subvariety of a bundle over the Teichm\"{u}ller space and is not connected as seen below. $\mathcal{M}^n_g(\textbf{Q})$ is a suitable subset of $\mathcal{M}^n_g$ and the orthogonal group $\mathcal{O}(n)$ acts on $\mathcal{M}^n_g(\textbf{Q})$. The Moduli space is defined as a quotient set $\mathcal{M}^n_g(\textbf{Q})/\mathcal{O}(n)$. In particular, for $n=4$, it is possible to give precise description of $\mathcal{M}^4_g$ \cite{A-Pi}. C. Arezzo and G. P. Pirola gave four components of $\mathcal{M}^4_g$ and calculated their dimensions. The details are as follows. A component corresponding to holomorphic minimal surfaces has complex dimension $5g-2$, a component projecting onto the hyperelliptic curves has dimension $4g$, a component projecting onto the non-hyperelliptic curves has dimension $4g$, and a component projecting onto general curves has dimension $4g$. They left the existence problem open for other components of $\mathcal{M}^4_g$. In this paper, we consider trigonal or $d$-gonal minimal surfaces in flat $4$-tori. A Riemann surface is said to be $d$-gonal if it can be represented as a $d$-sheeted branched cover of the sphere. The set of $d$-gonal Riemann surfaces is called $d$-gonal locus. In particular, we call trigonal for the case $d=3$. Then, our first main result is the following:
\begin{main}$\,$\\
(1) $\mathcal{M}^4_g$ has components projecting onto hyperelliptic locus distinct from Arezzo-Pirola's component. These components have complex $4g$-dimension.\\
(2) $\mathcal{M}^4_g$ has a component projecting onto trigonal locus. This component has a complex $4g$-dimension.\\
(3) Let $g>2(d-1)$. Then, $\mathcal{M}^4_g$ has a component projecting onto $d$-gonal locus and this component has a complex $4g$-dimension.
\end{main}
These results suggest that a component of $\mathcal{M}^4_g$ would have either dimension $4g$ or $5g-2$. Hence, the component corresponding to holomorphic minimal surfaces may give the maximal dimension.

Before we refer to (ii), we give some comments about associate surfaces. A compact minimal surface $M_g$ in a flat torus $\textbf{R}^n/\Lambda$ can be replaced by an $n$-periodic minimal surface in $\textbf{R}^n$ and the coodinate functions are harmonic. The classical result is the following (p.884 \cite{Me}):
\begin{thm}{\rm (Generalized  Weierstrass Representation)}
If $f:M_g\longrightarrow\textbf{R}^n/\Lambda$ is a conformal minimal immersion then, after a translation, $f$ can be represented by
\begin{equation*}
f(p)=\Re\int_{p_0}^p (\omega_1, \omega_2, \dots, \omega_n)^T \quad Mod \>\> \Lambda,
\end{equation*}
where $p_0 \in M_g$, superscript $T$ means the transposed matrix and $\{\omega_1$, $\omega_2$, $\dots$, $\omega_n \}$ are holomorphic differentials on $M_g$ satisfying
\begin{align}
&\{\omega_1, \omega_2, \dots, \omega_n\}{\rm\,\, has\,\, no\,\, common\,\, zeros},\\
&\sum_{i=1}^n \omega_i^2=0 \,,\\
&\left\{ \Re \int_{\gamma} (\omega_1, \omega_2, \dots, \omega_n)^T \,\,|\,\, \gamma \in H_1 (M_g, \textbf{Z}) \right\}{\rm is\,\, a\,\, sublattice\,\, of\,\, }\Lambda.
\end{align}
Conversely, every minimal surface in $\textbf{R}^n/\Lambda$ is obtained by the above construction.
\end{thm}
Condition (3) is called periodic condition and guarantees that the path integral is well-defined. Now we define the associate immersion $f_{\theta}$ of $f$:
\[f_{\theta}(p):=\Re\int_{p_0}^p e^{i\,\theta}\,(\omega_1, \omega_2, \dots, \omega_n)^T\]
If $f_{\theta}$ is well-defined, then we call it the associate surface of $f$. The conjugate surface of $f$ is the associate surface $f_{\pi/2}$ and its coodinate functions are the harmonic conjugates of the coordinate functions of $f$. The existence of the associate surface is an important subject for the theory of minimal surfaces. T. Nagano and B. Smyth \cite{N-S1} gave a criterion for the existence of the associate surfaces in the case $n=3$. In 1990, W. H. Meeks \cite{Me} introduced Property $\textbf{P}$ and considered the geometry of triply-periodic minimal surfaces that have an infinite number of distinct isometric minimal immersions in flat $3$-tori. We recall that a minimal surface $f:M\longrightarrow T^3$ is said to satisfy Property $\textbf{P}$ if for a countable dense set of angle $\theta\subset S^1$, the associate surface $f_{\theta}$ exists.

Now we explain about (ii). 
It is also an important subject to give an example of minimal surfaces. 
By the way, our objects are hyperelliptic, trigonal, and $d$-gonal minimal surfaces 
in flat $4$-tori. Thus, it is natural to ask if there exists an example ? 
We can construct hyperelliptic minimal surfaces by classical examples. 
Hence, we will deal with trigonal minimal surfaces. 
It is known that the genus $g$ of trigonal Riemann surfaces is $g\geq4$. 
We consider in the simplest case $g=4$. However, we first construct a minimal surface 
in a flat $3$-torus satisfying Property $\textbf{P}$ since Meeks proved that such 
a minimal surface can be full minimally immersed into 
a flat $4$-torus (Corollary 9.1 \cite{Me}). 
Now we review some notations and compare our example with classical examples. 
Let $f:M\longrightarrow T^3$ be a minimal surface in a flat $3$-torus. 
Recall that the space group $S_f (M)$ of $M$ is the set of symmetries of $M$ in $T^3$ 
(Definition 5.5 , Definition 6.1 \cite{Me}). 
We will denote the orientation-preserving subgroup of $S_f (M)$ by $S_f^o (M)$ and 
orientation-reversing elements by $S_f^r (M)$. 
Many classical examples of minimal surfaces in flat $3$-tori have cubical lattices for 
their tori, and have the standard linear symmetries of a cube in their space group. 
We will say that a minimal surface $f:M\longrightarrow \textbf{R}^3/\Lambda$ has 
the symmetries of a cube if $S_f^o (M)$ has members with linear part
\[R_1=\begin{pmatrix}0&0&1\\
                     1&0&0\\
                     0&1&0 \end{pmatrix}\]
and
\[ R_2=\begin{pmatrix} 1&0&0\\
                       0&-1&0\\
                       0&0&-1 \end{pmatrix}
\]
and if $S_f^r (M)$ has an inversion symmetry $T$ with linear part of $T$ equaling $-I_3$, 
where $I_3$ means the $3\times 3$ identity matrix.
Meeks proved the following:
\begin{thm}{\rm (}Theorem 6.2 \cite{Me}{\rm )}

Let $f:M\longrightarrow T^3$ be a minimal surface. Suppose that

 (1) $f$ can be lifted to a holomorphic immersion $f:M\longrightarrow \textbf{C}^3/\Lambda$,

(2) $f$ has the linear symmetries of a cube.

Then, $f$ satisfies Property $\textbf{P}$.
\end{thm}
Several classical examples of minimal surfaces satisfying Property $\textbf{P}$ 
satisfy the conditions of the theorem, e.g. 
the Schwarz surface \cite{Sch}, and A. Schoen's $O$, $C-TO$ surface \cite{A}. 
We construct a trigonal minimal surface having the conjugate surface, 
satisfying Property $\textbf{P}$. However, our trigonal minimal surface does not 
have the linear symmetries of a cube. Now we state our second main result:
\begin{main}$\,$\\
(1) There exists an example of trigonal minimal surfaces of genus $4$ in $3$-tori having the conjugate surface, satisfying Property $\textbf{P}$.\\
(2) The symmetry condition in Meeks' theorem is not always necessary for Property $\textbf{P}$.
\end{main}
As a corollary of Main Theorem 2, we otain the following:
\begin{cor}$\,$

There exits a trigonal minimal surface of genus $4$ in a flat $4$-torus. Therefore, the subset of $\mathcal{M}^4_4(\textbf{Q})/\mathcal{O}(4)$ given by the component in Theorem 1 (2) is not empty.
\end{cor}

Finally, we note that many important results on periodic minimal surfaces pass through 
the understanding of the space of Jacobi fields of a given one. It is an interesting problem 
to calculate the index and nullity of our trigonal minimal surface. We will set it aside 
and leave it as another problem.

The author would like to thank A. Futaki and the referee for their useful comments.

\section{The Moduli space of minimal surfaces in flat tori}

First of all, we review some fundamental results in \cite{A-Pi}. 
We then consider a fixed symplectic basis $\{\alpha_i,\beta_i\}_{i=1}^{g}$ of $H_1(M_g,\textbf{Z})$, and define
\[\mathcal{T}_g:=\{(M_g,\{\alpha_i,\beta_i\}_{i=1}^{g})\}/\sim,\]
where $(M_g,\{\alpha_i,\beta_i\}_{i=1}^{g})\sim(M'_g,\{\alpha_i,\beta_i\}_{i=1}^{g})$ if and only if there exists a biholomorphism $\phi:M_g \longrightarrow M'_g$ such that $\phi_{*}([\alpha_i])=[\alpha_i]$ and $\phi_{*}([\beta_i])=[\beta_i]$ ($[\>\>\> ]$ denotes the class in homotopy). We also define
\begin{align*}
\mathcal{H}^n_g:=\{(M_g,\{\alpha_i,\beta_i\}_{i=1}^{g}), \omega_1, \omega_2, \cdots, \omega_n\,|\,(M_g,\{\alpha_i,\beta_i\}_{i=1}^{g})\in\mathcal{T}_g, \\
\quad\omega_1,\omega_2, \cdots, \omega_n \in H^0(M_g, K)\,\, {\rm without} \,\,{\rm common}\,\,{\rm zeros}\},
\end{align*}
and
\begin{align*}
\mathcal{M}^n_g:=\{(M_g,\{\alpha_i,\beta_i\}_{i=1}^{g}), \omega_1, \omega_2, \cdots, \omega_n\in\mathcal{H}^n_g\,|\,\omega_1^2+\omega_2^2+ \cdots + \omega_n^2=0\}.
\end{align*}
Clearly $\mathcal{H}^n_g$ is a smooth manifold of complex dimension $3g-3+ng$. It is then natural to study the period mapping
\[\Pi((M_g,\{\alpha_i,\beta_i\}_{i=1}^{g}), \omega_1, \omega_2, \cdots, \omega_n):=\Re\scriptsize{\begin{pmatrix} \int_{\alpha_1}\omega_1 & \int_{\beta_1}\omega_1 & \cdots &\int_{\alpha_g}\omega_1 & \int_{\beta_g}\omega_1\\
\int_{\alpha_1}\omega_2 & \int_{\beta_1}\omega_2 & \cdots &\int_{\alpha_g}\omega_2 & \int_{\beta_g}\omega_2\\
\cdots&\cdots&\cdots&\cdots&\cdots\\
\int_{\alpha_1}\omega_n & \int_{\beta_1}\omega_n & \cdots &\int_{\alpha_g}\omega_n & \int_{\beta_g}\omega_n
\end{pmatrix}}\]
defined on $\mathcal{M}^n_g$. Observe that $\mathcal{T}_g$ can be identified via the classical period mapping with the Jacobian locus inside the Siegel Moduli space. We observe that condition $(3)$ if and only if
\begin{equation}
rank_{\textbf{Q}}\Pi((M_g,\{\alpha_i,\beta_i\}_{i=1}^{g}), \omega_1, \omega_2, \cdots, \omega_n)=n
\end{equation}

We can consider 
$\Pi((M_g,\{\alpha_i,\beta_i\}_{i=1}^{g}), \omega_1, \omega_2, \cdots, \omega_n)$ 
as the element of Grassmannian of $n$-planes in $H_1(M_g, \textbf{R})$. 
In fact, the Grassmannian $G(n, 2g)$ of $n$-planes in $\textbf{R}^{2g}$ can be represented 
by $M(n, 2g; n)/\sim$, where $M(n, 2g; n)$ denotes $(n, 2g)$ real matrices of rank $n$ and 
$A \sim B$ if and only if the $n$-plane spanned by $A$ agrees with $B$. 

$(4)$ means that there are $2g-n$ linear combinations of the columns of $\Pi((M_g,\{\alpha_i,\beta_i\}_{i=1}^{g})$, $\omega_1, \omega_2, \cdots, \omega_n)$ with rational coefficients. Thus, we define
\begin{align*}
&\mathcal{G}_{\textbf{Q}}, \,\,{\rm the\,\, Grassmannian\,\, of}\,\,n{\rm-planes\,\, in}\,\,H_1(M_g,\textbf{Q}),\\
&\mathcal{M}^n_g(\textbf{Q}):=\{p \in \mathcal{M}^n_g \,|\,\Pi(p) \in \mathcal{G}_{\textbf{Q}}\}.
\end{align*}
Observe that $\mathcal{M}^n_g(\textbf{Q})$ is the set of points in $\mathcal{M}^n_g$ for which a minimal immersion in a flat tori is well defined, or equivalently, is the set of points such that $\Pi(p)$ has rank $n$ over $\textbf{Q}$.

From the Riemannian geometric point of view it is natural to identify minimal immersions obtained by compositions with isometries of the target tori. This suggests defining an action of $\mathcal{O}(n)$ on $\mathcal{M}^n_g$ by
\[O((M_g,\{\alpha_i,\beta_i\}_{i=1}^{g}), \omega_1, \omega_2, \cdots, \omega_n)=((M_g,\{\alpha_i,\beta_i\}_{i=1}^{g}), (\omega_1, \omega_2, \cdots, \omega_n)O).\]
We will then call the space $\mathcal{M}^n_g(\textbf{Q})/\mathcal{O}(n)$ the Moduli space of minimal immersions of Riemann surfaces of genus $g$.

\section{proof of the Main Theorem 1}

In this section, we review the surveys of $\mathcal{M}^4_g$ in \cite{A-Pi} and give a proof of our main result. For any point $p$ in $\mathcal{M}^n_g$ we can consider the Gauss map $G:M_g\longrightarrow \textbf{C}P^{n-1}$ given by the evaluation of the holomorphic differentials. Because of equation $(2)$, we know that the image of the Gauss map is contained in the quadric $Q_{n-2}:=\{w \in \textbf{C}P^{n-1}|w \cdot w=0\}$, where "$\,\cdot\,$" is the complex bilinear inner product. If $n=4$, $Q_2$ is isomorphic to the product of two lines, and isomorphism is given by the Veronese map $V$ (p.19 \cite{H-Os}). We have the following commutative diagram:
$$
\begin{diagram}
M_g	& 		& \rTo^G	&		& \textbf{C}P^3 \\
	& \rdTo^{\varphi}	&		& \ruTo^{V}	&   \\
	&		& \textbf{C}P^1\times \textbf{C}P^1		&		&   \\
\end{diagram}
$$

$V$ can be given by
\[V((s_1,s_2), (t_1,t_2))=(s_1\,t_1-s_2\,t_2,\, s_1\,t_2+s_2\,t_1,\,i\,(s_1\,t_1+s_2\,t_2),\,i\,(s_1\,t_2-s_2\,t_1)).\]

Let us call $\varphi_1$ and $\varphi_2$ the projection of $\varphi$ on each of the two factors. Because the Gauss map is evaluation of holomorphic differentials, we know that $G^{*}(\mathcal{O}(1))=[(\omega_i)]=K$. If we put $\varphi_i^{*}(\mathcal{O}(1))=L_i$ ($i=1,\,2$), the commutativity of the diagram implies that $L_1+L_2=K$. This pair of line bundles still does not identify the minimal immersion, because we are not able to reconstruct the holomorphic differentials. This can be done precisely if we know not just the line bundles but also a pair of holomorphic sections for each of them. In fact, $\varphi_1=[s_1, s_2]$, and $\varphi_2=[t_1,t_2]$, where $s_i \in H^0(M_g, L_1)$ and $t_i \in H^0(M_g, L_2)$ ($i=1,\,2$). Knowing these sections, we can recover the holomorphic differentials just using the expression of the Veronese map.

We define
\begin{align*}
\widetilde{\mathcal{M}}^4_g:=\{(M_g, \{\alpha_i,\beta_i\}_{i=1}^{g}), L, &s_1, s_2, t_1, t_2\,|\\
&s_i \in H^0(M_g, L),\,t_i \in H^0(M_g, K-L)\},
\end{align*}
and construct a surjection $\Psi:\widetilde{\mathcal{M}}^4_g\longrightarrow \mathcal{M}^4_g$ given by
\begin{align*}
\omega_1&=s_1\,t_1-s_2\,t_2,\\
\omega_2&=s_1\,t_2+s_2\,t_1,\\
\omega_3&=i\,(s_1\,t_1+s_2\,t_2),\\
\omega_4&=i\,(s_1\,t_2-s_2\,t_1).
\end{align*}
A simple calculation shows that
\begin{align*}
\Psi^{-1}&((M_g, \{\alpha_i,\beta_i\}_{i=1}^{g}), \omega_1, \omega_2, \omega_3, \omega_4)\\
&=\{(M_g, \{\alpha_i,\beta_i\}_{i=1}^{g}), L, \lambda\,s_1, \lambda\,s_2, \lambda^{-1}\,t_1, \lambda^{-1}\,t_2\,|\,\lambda \in \textbf{C}^{*}\}
\end{align*}
and therefore
\[\mathcal{M}^4_g \cong \widetilde{\mathcal{M}}^4_g/\textbf{C}^{*}.\]

In this paper we restrict ourselves to the study of minimal immersions free of branch points. A direct calculation shows that this is equivalent to $s_1=s_2=0$ or $t_1=t_2=0$. $L$ and $K-L$ have two holomorphic sections without common zeros. This means $|D|$ and $|K-D|$ are base point free, where $D$ is the divisor coressponding to $L$ and $K$ means the canonical divisor. Therefore we can assume $h^0(L)>1$ and $h^0(K-L)>1$, unless one of these bundles is trivial. It is easy to see that this can happen if and only if the immersion is holomorphic with respect to a complex structure of the tori. 

 C. Arezzo and G. P. Pirola \cite{A-Pi} gave four components of $\mathcal{M}^4_g$. For example, a component corresponding to holomorphic minimal surfaces has complex dimension $5g-2$, a component projecting onto hyperelliptic locus and containing the spin bundles has dimension $4g$, a component projecting onto non-hyperelliptic locus and containing the spin bundles has dimension $4g$, and a component projecting onto general curves has dimension $4g$. Now we give some comments about the last component. 
We introduce the variety $W^r_d(M_g)$ (p.153, p.176 \cite{A-C-G-H}):
\begin{align*}
supp(W^r_d (M_g))&=\{L \in Pic^d (M_g)\,|\,h^0(L)\geq r+1\}\\
                 &=\{|D| \,|\, deg\, D=d\>,\>\>h^0(D)\geq r+1,\}
\end{align*}
where $r\geq d-g$.
\begin{lem} {\rm (p.182 \cite{A-C-G-H})}
Suppose $r\geq d-g$. Then no component of $W^r_d (M_g)$ is entirely contained in $W^{r+1}_d (M_g)$. In particular, if $W^r_d (M_g)\neq \phi$, then $W^r_d (M_g)-W^{r+1}_d (M_g)\neq \phi$.
\end{lem}
We consider $W^r_d (M_g)$ on general curves. If the Brill-Noether number $\rho=g-(r+1)\,(g-d+r)\geq 0$, then $W^r_d (M_g)\neq \phi$ and $W^r_d (M_g)-W^{r+1}_d (M_g)$ is smooth of dimension $\rho$ (p.190, p.214 \cite{A-C-G-H}). The last component is given by $deg\, L=d$, $h^0(L)=2$. Note that $L$ varies in $W^1_d (M_g)-W^{2}_d (M_g)$ in this case. Hence, we can show that the dimension is $4g$ (see p.775 \cite{A-Pi}).

They leave open the question of existence for other components of $\mathcal{M}^4_g$.

For general curves, if the Brill-Noether number $\rho=g-(r+1)\,(g-d+r)< 0$, then $W^r_d (M_g)$ is empty (p.214 \cite{A-C-G-H}). It is natural to ask whether $W^r_d (M_g)$ is empty for negative $\rho$ on other curves. In this generality the answer to this question is obviously no. Hyperelliptic curves in genus $g\geq3$, trigonal curves in genus $g\geq5$, and $d$-gonal curves in genus $g>2(d-1)$, give examples for which $\rho<0$ but $W^r_d (M_g)$ is empty (p.212 \cite{A-C-G-H}). We now consider components projecting onto these curves.\\
\\
(Components projecting onto the hyperelliptic locus)

We first give a notation about linear series (or system) \cite{A-C-G-H}. Let $V$ be a vector subspace of $H^0(M_g, \mathcal{O}(D))$ and $\mathcal{P}V$ is the projective space. A linear series $\mathcal{P}V$ is said to be a $g^r_d$ if $deg\, D=d$, $\dim V=r+1$.

Let $d$ be an integer such that $0\leq d\leq g$. Then any complete linear series $g^r_d$ on hyperelliptic Riemann surface $M_g$ is of the form
\[r\,g^1_2+p_1+p_2+\cdots+p_{d-2r},\]
where no two of the $p_i$'s are conjugate under the hyperelliptic involution (p.13 \cite{A-C-G-H}). We may assume $0\leq deg\, L\leq g-1$, if necessary we reconsider $L\longrightarrow K-L$. Since we now deal with the base point free cases, we obtain the form $g^r_d=r\,g^1_2$. In particular, $d=2r$. It follows that $W^r_{2r} (M_g)-W^{r+1}_{2r} (M_g)$ consists of the single point $g^r_{2r}$. The components of $\mathcal{M}^4_g$ given by $deg\,L=2r$, $h^0(L)=r+1$ then have dimension as follows:
\begin{align*}
2g-1+0+2h^0(L)+2h^0(K-L)-1&=2g-1+2(r+1)+2(g-r)-1\\
&=4g,
\end{align*}
where we use Riemann-Roch.
\begin{rem}
Arezzo-Pirola's component projecting hyperelliptic locus is in the case $r=\dfrac{g-1}{2}$ {\rm (}in the case that the genus is odd {\rm )}.
\end{rem}
$\,$\\
(Components projecting onto the trigonal locus)

Note that a Riemann surface $M_g$ is trigonal if and only if 
it has a base point free linear series $g^1_3$. 
Moreover, it is known that the genus of trigonal curve is $g\geq 4$ and 
the trigonal locus has complex dimension $2g+1$. We consider $W^1_3 (M_g)$ on trigonal curves. 
Note that the Brill-Noether number $\rho=g-2(g-3+1)=-g+4\leq 0$. 
The component of $\mathcal{M}^4_g$ given by $deg\,L=3$, $h^0(L)=2$ 
then has dimension as follows:

$\underline{{\rm the}\> {\rm case}\> {\rm of}\> g=4}$

It is known that $W^1_3 (M_g)$ consists of two distinct points or a single point 
(p.206 \cite{A-C-G-H}). By Lemma 3.1, $W^1_3 (M_g)-W^2_3 (M_g)$ consists of these points. 
Thus, we can calculate the dimension as follows:
\[(2 \times 4+1)+0+2\times2+2\times2-1=16=4g.\]
\begin{rem}
In the previous case {\rm (}hyperelliptic case{\rm)}, 
we note that we may consider in case $0\leq deg L\leq g-1$.  
If $M_g$ is neither hyperelliptic nor L=0, then we can reduce to in case 
$2r<deg L\leq g-1$ by Clliford's Theorem {\rm (}p.107 \cite{A-C-G-H}{\rm )}. 
In case $g=4$, the only values of $deg L$ and $r$ to be considered are $deg L=3$ and $r=1$. 
Therefore, components of $\mathcal{M}^4_g$ on trigonal lucus consists of $W^1_3 (M_g)$ alone 
in this case.
\end{rem}

$\underline{{\rm the}\> {\rm case}\> {\rm of}\> g>4}$

It is well known that a trigonal curve has unique complete $g^1_3$ for $g>4$ (p.37 \cite{C-M}, p.244 \cite{C-K-M}). Thus, $W^1_3 (M_g)-W^2_3 (M_g)$ consists of the $g^1_3$. Hence, we can calculate the dimension as follows:
\[(2g+1)+0+2\times2+2(g-2)-1=4g.\]
\\
(Components projecting onto the $d$-gonal locus)

We consider $W^1_d (M_g)$ on $d$-gonal curves for $g>2(d-1)$. Then the Brill-Noether number $\rho=g-2(g-d+1)=-g+2(d-1)<0$. A base point free $g^1_d$ on $M_g$ defines a $d$-sheeted ramified covering of $\textbf{C}P^1$ which, by Riemann-Hurwitz Theorem, has $2d+2g-2$ ramification points. Since, up to projectivities, we can always assume that the $d$-sheeted covering in question is ramified over the points $0,\, 1,\,\infty$, we conclude that the $d$-gonal locus has dimension $2d+2g-5$. It is known that $g^1_d$ is unique (\cite{A-C}). By the uniqueness and Lemma 3.1, $W^1_d (M_g)-W^2_d (M_g)$ consists of the $g^1_d$. The component of $\mathcal{M}^4_g$ given by $deg\,L=d$, $h^0(L)=2$ then has dimension as follows:
\[2d+2g-5+0+2\times2+2(g-d+1)-1=4g.\]
These calculations complete the proof of Main Theorem 1.

\section{proof of the Main Theorem 2}

In this section, we give an example of trigonal minimal surfaces in $4$-dimensional flat tori. It is known that the genus $g$ of trigonal Riemann surfaces is $g\geq4$. We consider in the simplest case $g=4$. However, we first construct a trigonal minimal surface satisfying Property $\textbf{P}$ in a flat $3$-torus since Meeks proved the following:
\begin{thm}{\rm (}Corollary 9.1 \cite{Me}{\rm )}
Let $f:M\longrightarrow T^3$ be a minimal surface satisfying Property $\textbf{P}$. Then, $M$ can be fully immersed as a minimal surface in a flat $4$-torus.
\end{thm}

We will consider in a simple case for trigonal Riemann surfaces (for general case, see \cite{Buch}).
Let $a_1, a_2, \cdots, a_6$ be distinct complex values and $M$ a trigonal curve of genus $4$ defined by \[w^3=(z-a_1)(z-a_2)\cdots(z-a_6).\]
We can compute explicitly a basis for $H^0 (M, K)$, where $K$ means the canonical bundle on $M$. First note that we have an automorphism $\psi:M\longrightarrow M$ of order $3$ given by $\psi(z,w)=(z,\omega\, w)$, where $\omega=\dfrac{-1+\sqrt{3}\,i}{2}$. The induced linear transformation \[\psi^*:H^0 (M, K)\longrightarrow H^0 (M, K)\]
is likewise of order $3$, and so a priori we obtain a decomposition of $H^0 (M, K)$ into eigenspaces with eigenvalues $1$, $\omega$, $\omega^2$. In fact, the $1$ eigenspace is trivial, since a holomorphic $1$-form $\eta$ on $M$ with $\psi^* \eta=\eta$ would descend to give a holomorphic $1$-form on $\textbf{C}P^1$, and none such exists. Thus we have $\psi^* \eta=\omega\, \eta$ or $\psi^* \eta=\omega^2\, \eta$ for all $\eta\in H^0 (M, K)$.

Now consider the $1$-form \[\omega_0=\dfrac{dz}{w^2}\] on $M$. $\omega_0$ is holomorphic and non-zero away from the points $z=\infty$, since the points where $w^2$ vanishes are exactly the zeros of $dz$. Note that $\psi^* \omega_0=\omega\,\omega_0$. If $\eta$ is any other holomorphic $1$-form on $M$, \[\eta=h\cdot \omega_0,\] where $h$ is a meromorphic function  on $M$, holomorphic away from $z=\infty$. In the case $\psi^* \eta=\omega\, \eta$, we have $\psi^* h=h$, i.e., $h$ is a function of $z$ alone, and necessarily a polynomial in $z$. Thus candidates are $\omega_0$, $z \,\omega_0$, $z^2 \,\omega_0$, $z^3\, \omega_0,\cdots$. But $z^3\, \omega_0,\cdots$ have poles at $z=\infty$. Therefore, $\dfrac{dz}{w^2}, \,z\,\dfrac{dz}{w^2},\,z^2\,\dfrac{dz}{w^2}\in H^0 (M, K)$. Note that they have no common zeros. In a similar way, we obtain $\dfrac{dz}{w}\in H^0 (M, K)$ in the case $\psi^* \eta=\omega^2\, \eta$. Thus, we can write out a basis for $H^0 (M, K)$:\[\left\{\dfrac{dz}{w^2},\,z\,\dfrac{dz}{w^2},\,z^2\, \dfrac{dz}{w^2},\,\dfrac{dz}{w}\right\}.\]

Now let $M$ be a trigonal curve defined by $w^3=z^6-1$ and we consider a minimal immersion as follows:
\begin{align*}
f:&\>M\longrightarrow\textbf{R}^3\\
  &\>\>p\>\>\longmapsto\Re \int_{p_0}^{p} \left(\dfrac{1-z^2}{w^2},\,\dfrac{i\,(1+z^2)}{w^2},\,\dfrac{2\,z}{w^2}\right)^T\,dz.
\end{align*}
Then, we can obtain the period matrix $\Omega$ as follows (see Appendix 1):
\begin{align*}
     \Omega&= \Re \begin{pmatrix}
   X&Y
\end{pmatrix},\\
&=\small{\begin{pmatrix}
0&\dfrac{3}{2}\,A&-\dfrac{3}{2}\,A&0&-\dfrac{3}{2}\,A&0&\dfrac{3}{2}\,A&-\dfrac{3}{2}\,A\\
0&\dfrac{3}{2}\,B&\dfrac{3}{2}\,B&0&\dfrac{3}{2}\,B&0&-\dfrac{3}{2}\,B&-\dfrac{3}{2}\,B\\
-\sqrt{3}\,C&\dfrac{\sqrt{3}}{2}\,C&\dfrac{\sqrt{3}}{2}\,C&-\sqrt{3}\,C&-\dfrac{\sqrt{3}}{2}\,C&\sqrt{3}\,C&-\dfrac{\sqrt{3}}{2}\,C&-\dfrac{\sqrt{3}}{2}\,C
\end{pmatrix}},
\end{align*}
where $A$, $B$, $C$, $X$, and $Y$ are
\begin{align*}
\nonumber A&=\int^2_1 \dfrac{1}{((t^2-1)^2\,(4-t^2))^{\frac{1}{3}}}\,dt,\\
B&=\int^{\sqrt{3}}_0 \dfrac{2^{\frac{1}{3}}}{(t^2\,(3-t^2)^2)^{\frac{1}{3}}\,\sqrt{1+t^2}}\,dt,\\
\nonumber C&=\int^{\sqrt{3}}_0 \dfrac{2}{\left(t^2\,(3-t^2)^2\right)^{\frac{1}{3}}\,
\sqrt{4-t^2}}\,dt,
\end{align*}
\begin{align*}
X&=\begin{pmatrix}
(-\omega^2+\omega)\,A & (1-\omega)\,A & (-1+\omega^2)\,A & (-\omega^2+\omega)\,A\\
(\omega^2-\omega)\,B  & (1-\omega)\,B & (1-\omega^2)\,B & (-\omega^2+\omega)\,B \\
i\,(-\omega^2+\omega)\,C & i\,(1-\omega)\,C & i\,(-1+\omega^2)\,C & i\,(-\omega^2+\omega)\,C
\end{pmatrix}\\
Y&=\begin{pmatrix}
 (-1+\omega^2)\,A & (-\omega^2+\omega)\,A & (1-\omega)\,A & (-1+\omega^2)\,A\\
(1-\omega^2)\,B & (-\omega^2+\omega)\,B & (-1+\omega)\,B & (-1+\omega^2)\,B\\
i\,(1-\omega^2)\,C & i\,(\omega^2-\omega)\,C & i\,(-1+\omega)\,C & i\,(1-\omega^2)\,C
\end{pmatrix}.
\end{align*}
We review the lattice transformation. A lattice $\Lambda$ in a real vector space $\textbf{R}^n$ is a discrete subgroup of maximal rank in $\textbf{R}^n$. A flat tori is a quotient $\textbf{R}^n/\Lambda$ with a lattice $\Lambda$ of $\textbf{R}^n$ and the metric induced from the standard Euclidean metric on $\textbf{R}^n$. Let $\{u_1,\cdots, u_m\}$ ($m \geq n$) be a sequence of vectors which span $\textbf{R}^n$. In general, $\{u_1,\cdots, u_m\}$ are not lattice vectors.
\begin{pro}{\rm (\cite{E} section 6)}
$\{u_1,\cdots, u_m\}$ are lattice vectors if and only if there exist lattice vectors $\{v_1,\cdots, v_n\}$ such that
\begin{align}
\nonumber(v_1, v_2, \cdots, v_n)=(u_1, u_2, \cdots, u_m)\,G_1,\\
\nonumber(u_1, u_2, \cdots, u_m)=(v_1, v_2, \cdots, v_n)\,G_2,
\end{align}
where $G_1$ is an $(m,n)$-matrix and $G_2$ is an $(n,m)$-matrix whose components are integers.
\end{pro}
Now we put the following three matrices:
\begin{align*}
\Lambda:=
\begin{pmatrix}
3A&0&\dfrac{3}{2}\,A\\
0&3B&\dfrac{3}{2}\,B\\
0&0&\dfrac{\sqrt{3}}{2}\,C
\end{pmatrix},
\>\>G_1:=\begin{pmatrix}
0&0&0\\
1&0&1\\
0&0&0\\
0&0&0\\
0&1&0\\
0&0&0\\
1&0&0\\
0&-1&0
\end{pmatrix},\\
G_2:=\begin{pmatrix}
1&0&-1&1&0&-1&1&0\\
1&0&0&1&1&-1&0&0\\
-2&1&1&-2&-1&2&-1&-1
\end{pmatrix},
\end{align*}
Then \[\Omega\,G_1=\Lambda,\>\>\Lambda\,G_2=\Omega.\]
By Proposition 4.1, $\Omega$ defines a lattice $\Lambda$ in $\textbf{R}^3$. 
Therefore, we can construct a trigonal minimal surface in a $3$-dimensional flat torus:
\begin{align}
f:&\>M\longrightarrow\textbf{R}^3/\Lambda\\
\nonumber  &\>\>p\>\>\longmapsto\Re \int_{p_0}^{p} \left(\dfrac{1-z^2}{w^2},\,\dfrac{i\,(1+z^2)}{w^2},\,\dfrac{2\,z}{w^2}\right)^T\,dz.\\
\nonumber&\qquad\qquad\qquad\qquad\qquad\qquad\qquad\qquad{\rm (}w^3=z^6-1{\rm )}
\end{align}
Next, we consider the conjugate surface and associate surface of $f$. 
By direct calculations, we can show that $f$ has the conjugate surface and 
satisfies Property $\textbf{P}$ (see Appendix 2). Hence, we obtain Main Theorem 2 (1).

Now, we prove Main Theorem 2 (2) as follows. 
First, $f$ satisfies $(1)$ in Theorem 1.2 because $f$ has the conjugate surface. 
However, $f$ does not satisfy $(2)$ in Theorem 1.2. 
We can see that by the reduction to absurdity. 

In fact, we assume that $(2)$ in Theorem 1.2 holds. 
Then $(3A,0,0)^T\in\Lambda$ implies $R_1(3A,0,0)^T=(0,3A,0)^T\in\Lambda$ and 
$R_1(0,3A,0)^T=(0,0,3A)^T\in\Lambda$. Since $(3A,0,0)^T$, $(0,3A,0)^T$, and 
$(0,0,3A)^T$ are linearly independent, $3A\,I_3$ defines the lattice in $\textbf{R}^3$ 
isomorphic to $\Lambda$, where $I_3$ is the $(3,3)$-identity matrix. Thus, by 
Proposition 4.1, there exist 
integral coefficient $(3,3)$-matrices $G_1$ and $G_2$ satisfying 
\[\Lambda=3A\, I_3\, G_1\>(=3A\,G_1),\>\>3A\, I_3=\Lambda\, G_2.\]
The first equation follows that  
\begin{align*}
G_1=\begin{pmatrix}
1&0&\dfrac{1}{2}\\
&*&\\
&**&\\
\end{pmatrix}, 
\end{align*}
and, therefore, it contradicts that $G_1$ consists of integers. 
As a result, we obtain Main Theorem 2.

Next our claim is Corollary 1.1. Applying Theorem 4.1 to $f$ defined by $(5)$, we can 
construct a full minimal surface in a $4$-dimensional flat torus. 
This minimal surface is not holomorphic with respect to any complex structure in the 
torus. This fact can be proved as follows:
By Meeks' construction, $f$ can be represented by
\[f(p)=\Re\int^p_{p_0} M\> \omega,\]
where $M$ and $\omega$ are defined by 
\[M=\begin{pmatrix}
1&0&0&0\\
0&1&0&0\\
0&0&b_1&0\\
0&0&b_2&0
\end{pmatrix}\>,\>\>\omega=
\begin{pmatrix}
\dfrac{1-z^2}{w^2}dz\\
\dfrac{i\,(1+z^2)}{w^2}dz\\
\dfrac{2\,z}{w^2}dz\\
\dfrac{dz}{w}
\end{pmatrix}.\]
On the other hands, Arezzo-Micallef \cite{A-Mi} showed necessary and sufficient 
condition for the holomorphy of minimal surfaces:
\begin{thm}
A full minimal immersion $f:M_g\longrightarrow \textbf{R}^n/\Lambda$ given by 
$f(p)=\Re\int^p_{p_0}M\omega$ is holomorphic with respect to some complex 
structure on the torus if and only if $M^TM=0$, where 
$\omega=(\omega_1,\cdots,\omega_g)^T$ is a basis of the holomorphic differential 
of $M_g$ and $M$ is a $(n,g)$ complex matrix.
\end{thm}
Now, it is easy to verify that $M^TM\neq0$ in our case. Hence, we obtain the 
non-holomorphy of $f$ by Theorem 4.2. 
Thus, the corresponding line bundle $L$ in $\mathcal{M}^4_g$ is neither $0$ nor $K$. 
Moreover, it is known that there exists no hyperelliptic minimal surface with even 
genus in a $3$-torus (Theorem 3.3 in \cite{Me}). 
It follows that our trigonal surface does not have the hyperelliptic 
structure. By Remark 3.2, we prove Corollary 1.1.

\section{Appendix 1}
We will calculate the period matrix in this section.
We introduce an automorphism $\varphi:M\longrightarrow M$ of order $3$ given by $\varphi(z,w)=(\omega z, w)$. We consider a basis for $1$-cycle on $M$ as follows:
\begin{align*}
A_1&=\left\{(z,\,w)=(e^{i\,t},\,w(t))\,|\,t\in \left[0,\dfrac{\pi}{3}\right],\,w\left(\dfrac{\pi}{6}\right)<0\right\}\\
&\cup \left\{(z,\,w)=(e^{-i\,t},\,\omega^2\,w(t))\,|\,t\in \left[-\dfrac{\pi}{3},\,0\right],\,w\left(-\dfrac{\pi}{6}\right)<0\right\}\\
A_2&=\varphi(A_1)=\left\{(z,\,w)=(e^{i\,t},\,w(t))\,|\,t\in \left[\dfrac{2\pi}{3},\,\pi\right],\,w\left(\dfrac{5\pi}{6}\right)<0\right\}\\
&\cup \left\{(z,\,w)=(e^{-i\,t},\,\omega^2\,w(t))\,|\,t\in \left[-\pi,\,-\dfrac{2\pi}{3}\right],\,w\left(-\dfrac{5\pi}{6}\right)<0\right\}
\end{align*}
\begin{align*}
A_3&=\psi(A_1)=\left\{(z,\,w)=(e^{i\,t},\,\omega\,w(t))\,|\,t\in \left[0,\dfrac{\pi}{3}\right],\,w\left(\dfrac{\pi}{6}\right)<0\right\}\\
&\cup \left\{(z,\,w)=(e^{-i\,t},\,w(t))\,|\,t\in \left[-\dfrac{\pi}{3},\,0\right],\,w\left(-\dfrac{\pi}{6}\right)<0\right\}\\
A_4&=\psi(A_2)=\left\{(z,\,w)=(e^{i\,t},\,\omega\,w(t))\,|\,t\in \left[\dfrac{2\pi}{3},\,\pi\right],\,w\left(\dfrac{5\pi}{6}\right)<0\right\}\\
&\cup \left\{(z,\,w)=(e^{-i\,t},\,w(t))\,|\,t\in \left[-\pi,\,-\dfrac{2\pi}{3}\right],\,w\left(-\dfrac{5\pi}{6}\right)<0\right\}\\
B_1&=\left\{(z,\,w)=(e^{i\,t},\,w(t))\,|\,t\in \left[-\pi,-\dfrac{2\pi}{3}\right],\,w\left(-\dfrac{5\pi}{6}\right)<0\right\}\\
&\cup \left\{(z,\,w)=(e^{-i\,t},\,\omega\,w(t))\,|\,t\in \left[\dfrac{2\pi}{3},\,\pi\right],\,w\left(\dfrac{5\pi}{6}\right)<0\right\}\\
B_2&=\varphi(B_1)=\left\{(z,\,w)=(e^{i\,t},\,w(t))\,|\,t\in \left[-\dfrac{\pi}{3},\,0\right],\,w\left(-\dfrac{\pi}{6}\right)<0\right\}\\
&\cup \left\{(z,\,w)=(e^{-i\,t},\,\omega\,w(t))\,|\,t\in \left[0,\,\dfrac{\pi}{3}\right],\,w\left(\dfrac{\pi}{6}\right)<0\right\}\\
B_3&=\psi(B_1)=\left\{(z,\,w)=(e^{i\,t},\,\omega\,w(t))\,|\,t\in \left[-\pi,-\dfrac{2\pi}{3}\right],\,w\left(-\dfrac{5\pi}{6}\right)<0\right\}\\
&\cup \left\{(z,\,w)=(e^{-i\,t},\,\omega^2\,w(t))\,|\,t\in \left[\dfrac{2\pi}{3},\,\pi\right],\,w\left(\dfrac{5\pi}{6}\right)<0\right\}\\
B_4&=\psi(B_2)=\left\{(z,\,w)=(e^{i\,t},\,\omega\,w(t))\,|\,t\in \left[-\dfrac{\pi}{3},\,0\right],\,w\left(-\dfrac{\pi}{6}\right)<0\right\}\\
&\cup \left\{(z,\,w)=(e^{-i\,t},\,\omega^2\,w(t))\,|\,t\in \left[0,\,\dfrac{\pi}{3}\right],\,w\left(\dfrac{\pi}{6}\right)<0\right\}
\end{align*}
We put \[\Phi=(\Phi_1,\Phi_2,\Phi_3)^T=\left(\dfrac{1-z^2}{w^2}\,dz,\,\dfrac{i\,(1+z^2)}{w^2}\,dz,\,\dfrac{2\,z}{w^2}\,dz\right)^T.\]
Then \[\psi^*\Phi=\omega\,\Phi,\>\>\varphi^*\Phi_3=\omega^2\,\Phi_3.\]
Hence
\begin{align*}
\int_{A_3}\Phi&=\omega\,\int_{A_1}\Phi,\,\int_{A_4}\Phi=\omega\,\int_{A_2}\Phi,\,\int_{B_3}\Phi=\omega\,\int_{B_1}\Phi,\\
\int_{B_4}\Phi&=\omega\,\int_{B_2}\Phi,\,\int_{A_2}\Phi_3=\omega^2\,\int_{A_1}\Phi_3,\,\int_{B_2}\Phi_3=\omega^2\,\int_{B_1}\Phi_3.
\end{align*}
Thus, it is sufficient to calculate \[\int_{A_1}\Phi,\>\> \int_{A_2}\Phi_1,\>\> \int_{A_2}\Phi_2, \>\>\int_{B_1}\Phi,\>\> \int_{B_2}\Phi_1,\>\> \int_{B_2}\Phi_2.\]
We will calculate the $A_1$-case. The other case can be solved in a similar way.\\
\\
\\
\\
$\underline{\Phi_1-case}$

We put $\eta=\dfrac{z^2+1}{2z}$. Then $d\eta=\dfrac{z^2-1}{2\,z^2}\,dz$ and
\[ \eta = \begin{cases}
     \dfrac{e^{2\,i\,t}+1}{2\,e^{i\,t}}=\cos t & (z=e^{i\,t})\\
     \cos t & (z=e^{-i\,t})
   \end{cases} \] Moreover
\[-\dfrac{2}{(4\,\eta^2-1)^2\,(1-\eta^2)}=8\,\dfrac{z^6}{(z^6-1)^2}=\left(2\,\left(\dfrac{z}{w}\right)^2\right)^3.\]
Now \[\dfrac{1-z^2}{w^2}\,dz=-2\,\left(\dfrac{z}{w}\right)^2\,\dfrac{z^2-1}{2\,z^2}\,dz.\]
We consider the case $t:0\longmapsto \dfrac{\pi}{3}$.
\begin{align*}
&2\,\left(\dfrac{z}{w}\right)^2\,\left(\dfrac{\pi}{6}\right)=\dfrac{2\,e^{\frac{\pi}{3}\,i}}{w\left(\frac{\pi}{6}\right)^2}=\dfrac{2\,(\omega+1)}{2^{\frac{2}{3}}}=-2^\frac{1}{3}\,\omega^2,\\
&-\dfrac{2}{(4\,\eta^2-1)^2\,(1-\eta^2)}\,\left(\dfrac{\pi}{6}\right)=-\dfrac{2}{(4\,\cos^2 \frac{\pi}{6}-1)^2\,(1-\cos^2 \frac{\pi}{6})}=-2.
\end{align*}
Hence, if we choose the branch as $\left(-\dfrac{2}{(4\,\eta^2-1)^2\,(1-\eta^2)}\right)^{\frac{1}{3}}<0$, then \[2\,\left(\dfrac{z}{w}\right)^2=\left(-\dfrac{2}{(4\,\eta^2-1)^2\,(1-\eta^2)}\right)^{\frac{1}{3}}\,\omega^2.\]
Thus
\[\dfrac{1-z^2}{w^2}\,dz=\left(\dfrac{2}{(4\,\eta^2-1)^2\,(1-\eta^2)}\right)^{\frac{1}{3}}\,\omega^2\,d\eta.\]
Note that $t:0\longmapsto \dfrac{\pi}{3}\Longleftrightarrow \eta:1\longmapsto \dfrac{1}{2}$. By evaluating $\eta=\dfrac{t}{2}$ ($t:2\longmapsto 1$), the integral of $\Phi_1$ is \[\int^1_2 \dfrac{2^{\frac{1}{3}}\,\omega^2}{((t^2-1)^2\,(1-\frac{t^2}{4}))^{\frac{1}{3}}}\,\dfrac{dt}{2}=\int^2_1 \dfrac{-\omega^2}{((t^2-1)^2\,(4-t^2))^{\frac{1}{3}}}\,dt.\]
Next, we consider the case $t:-\dfrac{\pi}{3}\longmapsto 0$. By the same argument as above, we obtain
\[\dfrac{1-z^2}{w^2}\,dz=\left(\dfrac{2}{(4\,\eta^2-1)^2\,(1-\eta^2)}\right)^{\frac{1}{3}}\,\omega\,d\eta.\]
Note that $t:-\dfrac{\pi}{3}\longmapsto 0\Longleftrightarrow \eta:\dfrac{1}{2}\longmapsto 1$. By evaluating $\eta=\dfrac{t}{2}$ ($t:1\longmapsto 2$), the integral of $\Phi_1$ is \[\int^2_1 \dfrac{2^{\frac{1}{3}}\,\omega}{((t^2-1)^2\,(1-\frac{t^2}{4}))^{\frac{1}{3}}}\,\dfrac{dt}{2}=\int^2_1 \dfrac{\omega}{((t^2-1)^2\,(4-t^2))^{\frac{1}{3}}}\,dt.\]
Adding the two integrals, we obtain
\[\int_{A_1}\dfrac{1-z^2}{w^2}\,dz=\int^2_1 \dfrac{-\omega^2+\omega}{((t^2-1)^2\,(4-t^2))^{\frac{1}{3}}}\,dt.\]
\\
$\underline{\Phi_2-case}$

We put $\eta=\dfrac{z^2-1}{z^2+1}$. Then $d\eta=\dfrac{4\,z}{(z^2+1)^2}\,dz$ and \[\eta=\begin{cases}
           \dfrac{e^{2\,i\,t}-1}{e^{2\,i\,t}+1}=i\,\tan t & (z=e^{i\,t})\\
            -i\,\tan t& (z=e^{-i\,t})
        \end{cases}\]
Moreover \[\left(\dfrac{z^2+1}{w}\right)^6=\left(\left(\dfrac{z^2+1}{w}\right)^2\right)^3=\left(\dfrac{4}{\eta^3+3\,\eta}\right)^2.\]
Now \[\dfrac{i\,(1+z^2)}{w^2}\,dz=i\,\left(\dfrac{z^2+1}{w}\right)^2\,\dfrac{z^2+1}{4\,z}\,\dfrac{4\,z}{(z^2+1)^2}\,dz.\]
We consider the case $t:0\longmapsto \dfrac{\pi}{3}$.
\begin{align*}
&\left(\dfrac{z^2+1}{w}\right)^2\,\left(\dfrac{\pi}{6}\right)=\dfrac{\left(e^{\frac{\pi}{3}\,i}+1\right)^2}{w\left(\frac{\pi}{6}\right)^2}=\dfrac{(\omega+2)^2}{2^{\frac{2}{3}}}=\dfrac{\omega^2+4\,\omega+4}{2^{\frac{2}{3}}}=-\dfrac{3\,\omega^2}{2^{\frac{2}{3}}},\\
&\left(\dfrac{4}{\eta\,(\eta^2+3)}\right)^2\,\left(\dfrac{\pi}{6}\right)=\dfrac{16}{-\tan^2 \frac{\pi}{6}\,(-\tan^2 \frac{\pi}{6}+3)^2}=-\dfrac{27}{4}.
\end{align*}
Hence, if we choose the branches as $\left(\dfrac{16}{\eta^2\,(\eta^2+3)^2}\right)^{\frac{1}{3}}=2\,\left(\dfrac{2}{\eta^2\,(\eta^2+3)^2}\right)^{\frac{1}{3}}<0$ and $\sqrt{1-\eta^2}>0$, then \[\left(\dfrac{z^2+1}{w}\right)^2=2\,\left(\dfrac{2}{\eta^2\,(\eta^2+3)^2}\right)^{\frac{1}{3}}\,\omega^2,\] and \[\dfrac{z^2+1}{4\,z}=\dfrac{e^{2\,i\,t}+1}{4\,e^{i\,t}}=\dfrac{\cos t}{2}=\dfrac{1}{2\,\sqrt{1-\eta^2}}\]
Thus \[\dfrac{i\,(1+z^2)}{w^2}\,dz=\dfrac{i\cdot2^{\frac{1}{3}}\,\omega^2}{\left(\eta^2\,(\eta^2+3)^2\right)^{\frac{1}{3}}\,\sqrt{1-\eta^2}}\,d\eta.\]
Note that $t:0\longmapsto \dfrac{\pi}{3}\Longleftrightarrow \eta:0\longmapsto \sqrt{3}\,i$. By evaluating $\eta=i\,t$ ($t:0\longmapsto \sqrt{3}$), the integral of $\Phi_2$ is \[\int^{\sqrt{3}}_0 \dfrac{i\cdot2^{\frac{1}{3}}\,\omega^2}{-(t^2\,(3-t^2)^2)^{\frac{1}{3}}\,\sqrt{1+t^2}}\,i\,dt=\int^{\sqrt{3}}_0 \dfrac{2^{\frac{1}{3}}\,\omega^2}{(t^2\,(3-t^2)^2)^{\frac{1}{3}}\,\sqrt{1+t^2}}\,dt.\]
Next, we consider the case $t:-\dfrac{\pi}{3}\longmapsto 0$. By the same argument as above, we obtain \[\dfrac{i\,(1+z^2)}{w^2}\,dz=\dfrac{i\cdot2^{\frac{1}{3}}\,\omega}{\left(\eta^2\,(\eta^2+3)^2\right)^{\frac{1}{3}}\,\sqrt{1-\eta^2}}\,d\eta.\]
Note that $t:-\dfrac{\pi}{3}\longmapsto 0\Longleftrightarrow \eta:\sqrt{3}\,i\longmapsto 0$. By evaluating $\eta=i\,t$ ($t:\sqrt{3}\longmapsto 0$), the integral of $\Phi_2$ is \[\int^{0}_{\sqrt{3}} \dfrac{i\cdot2^{\frac{1}{3}}\,\omega}{-(t^2\,(3-t^2)^2)^{\frac{1}{3}}\,\sqrt{1+t^2}}\,i\,dt=\int^{\sqrt{3}}_0 \dfrac{-2^{\frac{1}{3}}\,\omega}{(t^2\,(3-t^2)^2)^{\frac{1}{3}}\,\sqrt{1+t^2}}\,dt.\]
Adding the two integrals, we obtain \[\int_{A_1}\dfrac{i\,(1+z^2)}{w^2}\,dz=\int^{\sqrt{3}}_0 \dfrac{2^{\frac{1}{3}}\,(\omega^2-\omega)}{(t^2\,(3-t^2)^2)^{\frac{1}{3}}\,\sqrt{1+t^2}}\,dt.\]
\\
$\underline{\Phi_3-case}$

We put $\eta=\dfrac{z^2-1}{2\,z}$. Then $d\eta=\dfrac{z^2+1}{2\,z^2}\,dz$ and \[\eta=\begin{cases}
        \dfrac{e^{2\,i\,t}-1}{2\,e^{i\,t}}=i\,\sin t&z=e^{i\,t}\\
        -i\,\sin t&z=e^{-i\,t}
      \end{cases}\]
Moreover \[\left(\dfrac{z}{w}\right)^6=\left(\left(\dfrac{z}{w}\right)^2\right)^3=\dfrac{1}{64\,\eta^2\,(\eta^2+\frac{3}{4})^2}.\]
Now \[\dfrac{2\,z}{w^2}\,dz=2\,\left(\dfrac{z}{w}\right)^2\,\dfrac{2\,z}{z^2+1}\,\dfrac{z^2+1}{2\,z^2}\,dz.\]
We consider the case $t:0\longmapsto \dfrac{\pi}{3}$.
\begin{align*}
&\left(\dfrac{z}{w}\right)^2\,\left(\dfrac{\pi}{6}\right)=\dfrac{e^{\frac{\pi}{3}\,i}}{w\left(\frac{\pi}{6}\right)^2}=-\dfrac{\omega^2}{2^{\frac{2}{3}}},\\
&\dfrac{1}{64\,\eta^2\,(\eta^2+\frac{3}{4})^2}\left(\frac{\pi}{6}\right)=\dfrac{1}{-64\,\sin^2 \frac{\pi}{6}\,(-\sin^2 \frac{\pi}{6}+\frac{3}{4})^2}=-\dfrac{1}{4}
\end{align*}
Hence, if we choose the branches as $\left(\dfrac{1}{64\,\eta^2\,(\eta^2+\frac{3}{4})^2}\right)^{\frac{1}{3}}=\dfrac{1}{4\,\left(\eta^2\,(\eta^2+\frac{3}{4})^2\right)^{\frac{1}{3}}}<0$ and $\sqrt{1+\eta^2}>0$, then \[\left(\dfrac{z}{w}\right)^2=\dfrac{\omega^2}{4\,\left(\eta^2\,(\eta^2+\frac{3}{4})^2\right)^{\frac{1}{3}}},\] and \[\dfrac{2\,z}{z^2+1}=\dfrac{1}{\cos t}=\dfrac{1}{\sqrt{1+\eta^2}}.\]
Thus \[\dfrac{2\,z}{w^2}\,dz=\dfrac{\omega^2}{2\,\left(\eta^2\,(\eta^2+\frac{3}{4})^2\right)^{\frac{1}{3}}\,\sqrt{1+\eta^2}}\,d\eta.\]
Note that $t:0\longmapsto \dfrac{\pi}{3}\Longleftrightarrow \eta:0\longmapsto \dfrac{\sqrt{3}}{2}\,i$. By evaluating $\eta=\dfrac{i}{2}\,t$ ($t:0\longmapsto \sqrt{3}$), the integral of $\Phi_3$ is \[\int^{\sqrt{3}}_0 \dfrac{\omega^2}{-2\,\left(\frac{t^2}{4}\,(-\frac{t^2}{4}+\frac{3}{4})^2\right)^{\frac{1}{3}}\,\sqrt{1-\frac{t^2}{4}}}\,\dfrac{i}{2}\,dt=\int^{\sqrt{3}}_0 \dfrac{-2\,i\,\omega^2}{\left(t^2\,(3-t^2)^2\right)^{\frac{1}{3}}\,\sqrt{4-t^2}}\,dt.\]
Next, we consider the case $t:-\dfrac{\pi}{3}\longmapsto 0$. By the same argument as above, we obtain \[\dfrac{2\,z}{w^2}\,dz=\dfrac{\omega}{2\,\left(\eta^2\,(\eta^2+\frac{3}{4})^2\right)^{\frac{1}{3}}\,\sqrt{1+\eta^2}}\,d\eta.\]
Note that $t:-\dfrac{\pi}{3}\longmapsto 0\Longleftrightarrow \eta:\dfrac{\sqrt{3}}{2}\,i\longmapsto 0$. By evaluating $\eta=\dfrac{i}{2}\,t$ ($t:\sqrt{3}\longmapsto 0$), the integral of $\Phi_3$ is \[\int^{0}_{\sqrt{3}} \dfrac{\omega}{-2\,\left(\frac{t^2}{4}\,(-\frac{t^2}{4}+\frac{3}{4})^2\right)^{\frac{1}{3}}\,\sqrt{1-\frac{t^2}{4}}}\,\dfrac{i}{2}\,dt=\int^{\sqrt{3}}_0 \dfrac{2\,i\,\omega}{\left(t^2\,(3-t^2)^2\right)^{\frac{1}{3}}\,\sqrt{4-t^2}}\,dt.\]
Adding the two integrals, we obtain \[\int_{A_1}\dfrac{2\,z}{w^2}\,dz=\int^{\sqrt{3}}_0 \dfrac{2\,i\,(-\omega^2+\omega)}{\left(t^2\,(3-t^2)^2\right)^{\frac{1}{3}}\,\sqrt{4-t^2}}\,dt.\]

\section{Appendix 2}
In this section, we will show that $f$ has the conjugate immersion $f^*$, and moreover satisfies property $\textbf{P}$.
We consider the conjugate minimal immersion $f^*$ of $f$:
\begin{align*}
f^*:&\>M\longrightarrow\textbf{R}^3\\
  &\>\>p\>\>\longmapsto\Re \int_{p_0}^{p} i\,\left(\dfrac{1-z^2}{w^2},\,\dfrac{i\,(1+z^2)}{w^2},\,\dfrac{2\,z}{w^2}\right)^T\,dz.\\
&\qquad\qquad\qquad\qquad\qquad\qquad\qquad\qquad{\rm (}w^3=z^6-1{\rm )}
\end{align*}
Then, the period matrix $\Omega'$ is
\begin{align*}
     \Omega'&= \Re \begin{pmatrix}
   X'&Y'
\end{pmatrix},\\
&=\small{\begin{pmatrix}
-\sqrt{3}\,A&\dfrac{\sqrt{3}}{2}\,A&\dfrac{\sqrt{3}}{2}\,A&-\sqrt{3}\,A&\dfrac{\sqrt{3}}{2}\,A&-\sqrt{3}\,A&\dfrac{\sqrt{3}}{2}\,A&\dfrac{\sqrt{3}}{2}\,A\\
\sqrt{3}\,B&\dfrac{\sqrt{3}}{2}\,B&-\dfrac{\sqrt{3}}{2}\,B&-\sqrt{3}\,B&-\dfrac{\sqrt{3}}{2}\,B&-\sqrt{3}\,B&-\dfrac{\sqrt{3}}{2}\,B&\dfrac{\sqrt{3}}{2}\,B\\
0&-\dfrac{3}{2}\,C&\dfrac{3}{2}\,C&0&-\dfrac{3}{2}\,C&0&\dfrac{3}{2}\,C&-\dfrac{3}{2}\,C
\end{pmatrix}},
\end{align*}
where $X'$ and $Y'$ are
\begin{align*}
X'&=\begin{pmatrix}
i\,(-\omega^2+\omega)\,A & i\,(1-\omega)\,A & i\,(-1+\omega^2)\,A & i\,(-\omega^2+\omega)\,A\\
i\,(\omega^2-\omega)\,B  & i\,(1-\omega)\,B & i\,(1-\omega^2)\,B & i\,(-\omega^2+\omega)\,B \\
(\omega^2-\omega)\,C & (-1+\omega)\,C & (1-\omega^2)\,C & (\omega^2-\omega)\,C 
\end{pmatrix}\\
Y'&=\begin{pmatrix}
 i\,(-1+\omega^2)\,A & i\,(-\omega^2+\omega)\,A & i\,(1-\omega)\,A & i\,(-1+\omega^2)\,A\\
i\,(1-\omega^2)\,B & i\,(-\omega^2+\omega)\,B & i\,(-1+\omega)\,B & i\,(-1+\omega^2)\,B\\
(-1+\omega^2)\,C & (-\omega^2+\omega)\,C & (1-\omega)\,C & (-1+\omega^2)\,C 
\end{pmatrix}.
\end{align*}
We put the following three matrices:
\begin{align*}
\Lambda':=
\begin{pmatrix}
\sqrt{3}A&0&\dfrac{\sqrt{3}}{2}\,A\\
0&\sqrt{3}B&\dfrac{\sqrt{3}}{2}\,B\\
0&0&\dfrac{3}{2}\,C
\end{pmatrix},
\>\>G'_1:=\begin{pmatrix}
0&1&0\\
1&1&-1\\
1&1&0\\
0&0&-1\\
0&0&0\\
0&0&0\\
0&0&0\\
0&0&0
\end{pmatrix},\\
G'_2:=\begin{pmatrix}
-1&1&0&-1&1&-1&0&1\\
1&1&-1&-1&0&-1&-1&1\\
0&-1&1&0&-1&0&1&-1
\end{pmatrix}.
\end{align*}
Then
\[\Omega'\,G'_1=\Lambda',\>\>\Lambda'\,G'_2=\Omega'\]
By Propositon 4.1, we can construct the conjugate minimal surface:
\begin{align*}
f^*:&\>M\longrightarrow\textbf{R}^3/\Lambda'\\
  &\>\>p\>\>\longmapsto\Re \int_{p_0}^{p} i\,\left(\dfrac{1-z^2}{w^2},\,\dfrac{i\,(1+z^2)}{w^2},\,\dfrac{2\,z}{w^2}\right)^T\,dz.\\
&\qquad\qquad\qquad\qquad\qquad\qquad\qquad\qquad{\rm (}w^3=z^6-1{\rm )}
\end{align*}

Next, we consider the associate surface $f_{\theta}$ of $f$:
\begin{align*}
f_{\theta}:&\>M\longrightarrow\textbf{R}^3\\
  &\>\>p\>\>\longmapsto\Re \int_{p_0}^{p} e^{i\,\theta}\,\left(\dfrac{1-z^2}{w^2},\,\dfrac{i\,(1+z^2)}{w^2},\,\dfrac{2\,z}{w^2}\right)^T\,dz.\\
&\qquad\qquad\qquad\qquad\qquad\qquad\qquad\qquad{\rm (}w^3=z^6-1{\rm )}
\end{align*}
Note that
\[\Re \int_{\gamma} e^{i\,\theta}\, \Phi=\cos \theta \,\Re \int_{\gamma}\Phi+\sin \theta \,\Re \int_{\gamma} i\,\Phi.\]
Hence, the period matrix $\Omega_{\theta}$ is given by
\begin{align*}
\Omega_{\theta}&=\cos \theta \,\Omega+\sin \theta\,\Omega'\\
&=\begin{pmatrix}
X_1&Y_1&Z_1
\end{pmatrix},
\end{align*}
where $X_1$, $Y_1$, and $Z_1$ are given by
\begin{align*}
X_1&=\small{\begin{pmatrix}
-\sqrt{3}\,\sin\theta\,A&\left(\dfrac{3}{2}\,\cos\theta+\dfrac{\sqrt{3}}{2}\,\sin\theta\right)A&\left(-\dfrac{3}{2}\,\cos\theta+\dfrac{\sqrt{3}}{2}\,\sin\theta\right)A\\
\sqrt{3}\,\sin\theta\,B&\left(\dfrac{3}{2}\,\cos\theta+\dfrac{\sqrt{3}}{2}\,\sin\theta\right)B&\left(\dfrac{3}{2}\,\cos\theta-\dfrac{\sqrt{3}}{2}\,\sin\theta\right)B\\
-\sqrt{3}\,\cos\theta\,C&\left(\dfrac{\sqrt{3}}{2}\,\cos\theta-\dfrac{3}{2}\,\sin\theta\right)C&\left(\dfrac{\sqrt{3}}{2}\,\cos\theta+\dfrac{3}{2}\,\sin\theta\right)C
\end{pmatrix}},\\
Y_1&=\small{\begin{pmatrix}
-\sqrt{3}\,\sin\theta\,A&\left(-\dfrac{3}{2}\,\cos\theta+\dfrac{\sqrt{3}}{2}\,\sin\theta\right)A&-\sqrt{3}\,\sin\theta\,A\\
-\sqrt{3}\,\sin\theta\,B&\left(\dfrac{3}{2}\,\cos\theta-\dfrac{\sqrt{3}}{2}\,\sin\theta\right)B&-\sqrt{3}\,\sin\theta\,B\\
-\sqrt{3}\,\cos\theta\,C&\left(-\dfrac{\sqrt{3}}{2}\,\cos\theta-\dfrac{3}{2}\,\sin\theta\right)C&\sqrt{3}\,\cos\theta\,C
\end{pmatrix}},\\
Z_1&=\small{\begin{pmatrix}
\left(\dfrac{3}{2}\,\cos\theta+\dfrac{\sqrt{3}}{2}\,\sin\theta\right)A&\left(-\dfrac{3}{2}\,\cos\theta+\dfrac{\sqrt{3}}{2}\,\sin\theta\right)A\\
\left(-\dfrac{3}{2}\,\cos\theta-\dfrac{\sqrt{3}}{2}\,\sin\theta\right)B&\left(-\dfrac{3}{2}\,\cos\theta+\dfrac{\sqrt{3}}{2}\,\sin\theta\right)B\\
\left(-\dfrac{\sqrt{3}}{2}\,\cos\theta+\dfrac{3}{2}\,\sin\theta\right)C&\left(-\dfrac{\sqrt{3}}{2}\,\cos\theta-\dfrac{3}{2}\,\sin\theta\right)C
\end{pmatrix}}.
\end{align*}
Now we choose $\theta$ as follows:
\begin{align}
2\sqrt{3}\sin\theta\,m&=(3\cos\theta+\sqrt{3}\sin\theta)\,n\\
\nonumber&\qquad\qquad\qquad\qquad{\rm (i.e.}\tan\theta=\dfrac{\sqrt{3}}{\frac{2\,m}{n}-1}\> {\rm )}
\end{align}
We may assume $n$ and $m$ are prime. There exist integers $x$ and $y$ such that $nx+my=1$. We now put
\begin{align*}
{\Omega'}_{\theta}&=\begin{pmatrix}
{{\Omega'}_{\theta}}_{1}&{{\Omega'}_{\theta}}_{2}
\end{pmatrix},\>\>\>
F_1=\begin{pmatrix}
-x&0&0&0&m&1&-1&0\\
y&0&0&0&n&1&0&0\\
0&0&0&0&0&1&0&0\\
0&0&0&0&0&0&1&0\\
0&-x&0&0&0&0&-1&m\\
-x&-x-y&1&0&m&1&-1&m-n\\
y&-x-y&0&1&n&1&-1&m-n\\
0&-y&0&0&0&1&0&-n\end{pmatrix},\\
F_2&=\begin{pmatrix}
-n&m&n-m&-n&0&0&0&0\\
0&0&m&-n&-n&0&0&-m\\
-1&0&1&-1&-1&1&0&-1\\
0&-1&1&0&-1&0&1&-1\\
0&0&x&y&y&0&0&-x\\
0&0&1&0&0&0&0&0\\
0&0&0&1&0&0&0&0\\
y&x&-x-y&y&0&0&0&0
\end{pmatrix}
\end{align*}
where  ${{\Omega'}_{\theta}}_{1}$, ${{\Omega'}_{\theta}}_{2}$ are given by
\begin{align*}
{{\Omega'}_{\theta}}_{1}&=\footnotesize{\begin{pmatrix}\dfrac{1}{n}\cdot2\sqrt{3}\sin\theta\,A&0&-\sqrt{3}\sin\theta\,A&\dfrac{m}{n}\cdot\sqrt{3}\sin\theta\,A\\
0&\dfrac{1}{n}\cdot2\sqrt{3}\sin\theta\,B&-\sqrt{3}\sin\theta\,B&-\dfrac{m}{n}\cdot\sqrt{3}\sin\theta\,B\\
0&0&\sqrt{3}\cos\theta\,C&\left(-\dfrac{\sqrt{3}}{2}\cos\theta+\dfrac{3}{2}\sin\theta\right)C
\end{pmatrix}}\\
{{\Omega'}_{\theta}}_{2}&=\begin{pmatrix}
0&0&0&0\\
0&0&0&0\\
0&0&0&0
\end{pmatrix}.
\end{align*}
Then \[\Omega_{\theta}\,F_1={\Omega'}_{\theta},\>\>{\Omega'}_{\theta}\,F_2=\Omega_{\theta}\]
By $(6)$,
\begin{align*}
{{\Omega'}_{\theta}}_{1}&=\footnotesize{\begin{pmatrix}\dfrac{1}{n}\cdot2\sqrt{3}\sin\theta\,A&0&-\sqrt{3}\sin\theta\,A&\dfrac{m}{n}\cdot\sqrt{3}\sin\theta\,A\\
0&\dfrac{1}{n}\cdot2\sqrt{3}\sin\theta\,B&-\sqrt{3}\sin\theta\,B&-\dfrac{m}{n}\cdot\sqrt{3}\sin\theta\,B\\
0&0&\sqrt{3}\cos\theta\,C&\left(-\dfrac{\sqrt{3}}{2}\cos\theta+\dfrac{3}{2}\sin\theta\right)C
\end{pmatrix}}\\
&=\footnotesize{\begin{pmatrix}\dfrac{1}{n}\cdot2\sqrt{3}\sin\theta\,A&0&-\sqrt{3}\sin\theta\,A&\dfrac{m}{n}\cdot\sqrt{3}\sin\theta\,A\\
0&\dfrac{1}{n}\cdot2\sqrt{3}\sin\theta\,B&-\sqrt{3}\sin\theta\,B&-\dfrac{m}{n}\cdot\sqrt{3}\sin\theta\,B\\
0&0&\sqrt{3}\cos\theta\,C&\dfrac{-m+2n}{2m-n}\cdot\sqrt{3}\cos\theta\,C
\end{pmatrix}}
\end{align*}
Thus, \[rank_{\textbf{Q}} \Omega_{\theta}=rank_{\textbf{Q}} {\Omega'}_{\theta}=3\]
It follows that $\Omega_{\theta}$ defines a lattice $\Lambda_{\theta}$ by $(4)$ and Proposition 4.1. Hence, we can also construct the associate surface $f_{\theta}$:
\begin{align*}
f:&\>M\longrightarrow\textbf{R}^3/\Lambda_{\theta}\\
  &\>\>p\>\>\longmapsto\Re \int_{p_0}^{p} e^{i\,\theta}\,\left(\dfrac{1-z^2}{w^2},\,\dfrac{i\,(1+z^2)}{w^2},\,\dfrac{2\,z}{w^2}\right)^T\,dz.\\
&\qquad\qquad\qquad\qquad\qquad\qquad\qquad\qquad{\rm (}w^3=z^6-1{\rm )}
\end{align*}
Moreover, by $(6)$, $f$ satisfies Property $\textbf{P}$.

\end{document}